\documentclass[12pt]{amsart}
\usepackage{amsthm,amsmath}

\def\tenv#1{\subsection{\bf #1} }
\newcommand{\ts}{\vspace{\baselineskip}\noindent{\bf Proof.}$\;\;$}

\def\n{{\noindent}}
\def\q{{\bf Q}}
\def\QQ{{\bf Q}}
\def\Q{{\q(\varphi)}}
\def\w{{\wedge}}
\def\wk{{\wedge_K}}
\def\L{{\bf Lemma}}
\def\CC{{\bf C}}
\def\f{{\varphi}}
\def\H{{\widetilde {H}}}
\def\T{{\bf Theorem}}
\def\D{{\bf Definition }}
\def\P{{\bf Proposition}}
\def\C{{\bf Corollary}}
\def\ep{{\epsilon}}

\def\io{{\iota}}
\def\RR{{\bf R}}
\def\ZZ{{\bf Z}}

\begin{document}
\title{ABELIAN VARIETIES OF WEIL TYPE AND KUGA-SATAKE VARIETIES}
\author{\sc Giuseppe Lombardo}
\date{\empty}

\maketitle
\begin{abstract} We analyze the relationship between abelian fourfolds of Weil type and Hodge structures of type K3, and we extend some of these correspondences to the case of  arbitrary dimension.\end{abstract}

\begin{center}
{\sc Introduction}
\end{center}

Abelian varieties of Weil type are examples of abelian
varieties for which the Hodge conjecture is open in general.
We study the structure of abelian fourfolds of Weil type, giving 
an explicit description of the Hodge structures of the 
cohomology groups, in particular the sub-Hodge structures of the 
second cohomology group.
Starting from the observations of Paranjape \cite{Pa}, we show that 
for certain abelian varieties of Weil type (those with discriminant one)  
there exists a polarized sub-Hodge structure of the second cohomology 
group 
of dimension 6 with $h^{2,0}=1$. We show that the map which 
associates this polarized Hodge structure to an abelian fourfold of Weil type with discriminant one  admits an 
``inverse''. Indeed, we can construct the Kuga-Satake variety 
associated to this Hodge structure and prove that the Kuga-Satake 
variety is an abelian variety of dimension 16 which is isogeneous to the product of 
four copies of 
the abelian fourfold of Weil type we started with. In the last 
section we generalize some of these results to higher dimensions;
starting from a polarized weight two Hodge structure of type 
$(1,n-2,1)$ with $n \equiv 2 \pmod 4 $, 
the Poincar\'e decomposition of the Kuga-Satake variety gives (a 
number of copies of) an abelian variety of Weil type with 
discriminant one.

I would like to thank my advisor, Professor B. van Geemen, for his continuing guidance and support
in the realization of this paper, and I thank the referee for his helpful comments.

\section{Abelian varieties of Weil type}
\tenv{Weil type} Let $(X,E)$ be a polarized abelian variety of 
dimension $2n$ and $K\hookrightarrow {\rm End}(X) \otimes \q$ be an imaginary quadratic 
field. The polarization $E\subset H^2(X,\q)$ defines by duality a \, bilinear antisymmetric map, which is for convenience denoted by the same letter $E$, from $H_1(X,\q) \times H_1(X,\q)$ to $\q$. $X$ is said to be of Weil 
type if the action of $K$ on the tangent space $T_{0}X$ can be diagonalized as
$${\rm diag} (\sigma (k), \dots ,\sigma (k), {\overline {\sigma(k)}}, 
\dots, {\overline {\sigma (k)}}) \qquad  (k \in K)$$
with $n$ entries $\sigma (k)$ and $n$ entries ${\overline {\sigma 
(k)}}$
(where $\sigma :K\hookrightarrow \CC$ is an embedding) and 

\n $E(kx,ky)= \sigma (k) \overline{ \sigma (k)}E(x,y)$ for $x,y \in T_0X$.

\tenv{Discriminant} Let$(X,K,E)$ be an abelian variety of Weil type 
and let $K=\Q$. The map 
$$\begin{array}{lrclcc}
H:& H_1(X,\q)& \times & H_1(X,\q) & \longrightarrow & K \\
 & (x& ,&y)&\longmapsto & E(\f x, y) +\f E(x,y)
\end{array}$$
is a Hermitian form on the $K$-vector space $H_{1}(X,\q)$. There 
exists a $K$-basis in which $H$ is represented by a diagonal matrix ${\rm diag} 
(a,1,\dots,1,-1,\dots, -1)$, where  $a$ is a rational positive number 
called the discriminant of the variety 
($a={\rm discr}(X,K,E)=(-1)^n {\rm det}(H)\in \q ^* / {\rm Nm}_{\q / K} (K^*)$).

\section{Hodge structures}
\tenv{} \label{Def}  Let $V$ be a $\q$-vector space and $h:S({\bf 
R}) \rightarrow
GL(V_{\bf R}) $ be a rational 
representation of the 
group $$S({\bf R})=\Big\{ s(a,b):=\left(\begin{array}{ll}a&b \\ -b&a 
\end{array}\right) 
\in GL(2,{\bf R}) \Big \} \cong \CC^{*}$$ on $V_{\bf R}=V \otimes_{\q} 
{\bf R}$ such that  $h(s(a,0))=a^n 
\cdot {\bf 1}$. The couple
$(V,h)$ is said to be a (rational) Hodge structure of weight $n$.

\tenv{}  By the action of $h$, we have a decomposition of $V_{\CC}=V 
\otimes_{\q}\CC$ into 
weight spaces:

$$V_{\CC}=
\bigoplus_{p+q=n}V^{p,q},$$

\n  where $V^{p,q}=\{ v \in V_{\CC} \ ;\  h(z)v=z^p\bar z^qv \}$ and 
$\overline{V^{q,p}}=V^{p,q}$. This decomposition of $V$ is a Hodge decomposition in the usual 
sense.

\tenv{Polarization}  A polarization of the Hodge structure $(V,h)$ is 
a  
$\q$-bilinear  map 

\n $\psi: V \times V \rightarrow \q$ such that

$\begin{array}{ll}
(i)&\psi (h(z)v,h(z)w)=(z \bar z)^n \psi (v,w) \quad {\rm for \ any}\  v,w \in 
V_{\bf 
R},\ z \in \CC^{*} \ \quad (n={\rm weight}),\\
(ii) & \psi (v,h(i)w)=\psi (w,h(i)v)\quad {\rm for \ any}\  v,w \in V_{\bf R},  \\
(iii) & \psi (v,h(i)v)>0 \quad {\rm for \ any}\  v \in V_{\bf R}-\{0 \}.
\end{array}$

\medskip

It is easy to show the following by direct computation (see \cite{vG}):

\tenv{\L} \label{pol} {\it A polarization $\psi$ is symmetric if the 
weight of the Hodge 
structure is even, alternating if the weight is odd, and the quadratic 
form $Q(v)=- \psi (v,v) $
associated to the polarization is positive definite on 
$(V^{2,0} \oplus V^{0,2}) \cap V_{\RR}$ and negative on $V^{1,1}\cap 
V_{\RR}$.}

\tenv{Weight 2 Hodge structures} Let $(V,h, \psi )$ be a polarized 
Hodge 
structure of weight $2$. The triple 
$({\rm dim}V^{2,0},{\rm dim}V^{1,1},{\rm dim}V^{0,2})$ 
is said to be the type of the Hodge structure. 

\section{Sub-Hodge structures of the cohomology groups}

\n Let $K := \Q$ ($\f^2=-d$) and $V$ be a $K$-vector space. We can extend the 
action of $K^*$ to $\bigwedge_{\q}^2 V$  in a natural way: $k_{2}(v 
\w w)=kv 
\w kw$ ($k \in K^{*}$). In order to examine the sub-Hodge structures, 
the following proposition concerning $K$-vector spaces is very useful (see also \cite{km}, \cite {mz} and \cite{cime} in which this proposition is implicitly contained):

\tenv{\P} {\label {naturale}}{\it Let $V$ be a $K$-vector space. Then the $\q$-linear map
$$
i:\bigwedge_{K}^2 V  \hookrightarrow  \bigwedge_{\q}^2 V, \qquad
a \w_{K} b  \mapsto  {1 \over 2} a \w b -{1 \over {2d}}\f a \w \f 
b
$$
is injective and ${\rm Im}(i)=\{w \in \bigwedge_{\q}^2 V \ ; \ \f_{2}
w=-dw\}$ {\rm (}where $\f_2$ is the action of $\f\in K^*$ on $\bigwedge^2_{\q}V$ already defined{\rm )}.}

\ts Let $W:={\rm Hom}_{K}(V,K)$, $W^* :={\rm Hom}_{K}(W,K)$ 
and $(W^*)^{\q }:={\rm Hom}_{\q}(W,\q)$.

\n The trace map 
$$\begin{array}{ccccccccccc}
\omega:&W^* & \stackrel {\cong}{\longrightarrow} &W^{*\q},&\quad &{\rm 
{where}}&\quad &{\rm Tr}:&K \ni z
&\longmapsto &z + \bar z \in \q\\ 
 &f &\longmapsto &{\rm Tr} \circ f&\quad & &\quad & &  
& &
\end{array}$$

\n is a linear isomorphism of $\q$-vector spaces.
We also have an isomorphism of $K$-vector spaces
$$\begin{array}{cccc}
\alpha:&\bigwedge^2_{K}V & \longrightarrow & 
{\rm Hom}_{K}(\bigwedge^2_{K}W, K) 
\\
 &a \w_{K} b & \longmapsto &\lbrack \alpha \w_{K} \beta \mapsto 
\alpha (a)\beta (b)-\beta (a)\alpha (b)
\rbrack. 
\end{array}$$ 

\n Then, since $K$-linear homomorphisms are also $\q$-linear, we 
have a chain of maps

$$\begin{array}{cccccccc} 
 &\bigwedge^2_{K}V &\stackrel{\alpha}\rightarrow & {\rm Hom}_{K} 
\Bigg(\bigwedge^2_{K}W,K \Bigg)
&\stackrel{\rm id}{\rightarrow}&{\rm Hom}_{\q} \Bigg( \bigwedge^2_{\q}W,K 
\Bigg) &\stackrel{\rm Tr}{\rightarrow}
&  \\
\stackrel{\rm Tr}{\rightarrow} &{\rm Hom}_{\q}(\bigwedge^2_{\q}W,\q) 
&\stackrel{\cong}\rightarrow 
&\bigwedge^2_{\q}(W^*)^{\q}&\stackrel{(\omega \w \omega)^{-1}} 
\longrightarrow&
\bigwedge^2_{\q}W^* &\stackrel{\cong}\rightarrow &\bigwedge^2_{\q}V.  
\end {array}$$  

\n Writing these maps in term of $\q$ and $K$-basis, we have that 
$i(a \w_{K} b)=(1/2) a \w b -(1/2d)\f a \w \f 
b.$ 

\n The map $i$ is a composite of injective maps, hence is 
injective.

Obviously,  ${\rm Im}(i) \subset \{w \in \bigwedge_{\q}^2 V \ ; \ \f_{2}
w=-dw\}$. Indeed, we have 
$$\f_{2}({1 \over 2} a \w b -{1 \over {2d}}\f a \w \f 
b)={1 \over 2} \f (a) \w \f (b) -{1 \over {2d}}\f^{2} (a) \w \f^{2} 
(b)=(-d)\bigg( {1 \over 2} a \w b -{1 \over {2d}}\f a \w \f b \bigg) 
.$$ 

\n Let now $V_{K}=V_{+} \oplus V_{-}$ be the decomposition of 
$V_{K}=V 
\otimes_{\q} K$ into the subspaces in which the elements $k \in K$ act 
respectively 
as $\sigma (k) $ and $\overline{ \sigma (k)}$; we have 
${\rm dim}_{K}V_{K}=2n$ and ${\rm dim}_{K}V_{\pm}=n$. By the action of  $K$ on  
$\bigwedge^{2}_{\q}V$, 
we have $\f_{2}^{2}=d^{2}$ and can decompose the space as 
$\bigwedge^{2}_{\q}V=W_{+} \oplus W_{-}$, where $W_{\pm}$ are the $+d$ 
and the $-d$ eigenspaces of $\f_2$. We showed that $${\rm Im}(i) \subset 
W_{-}$$
\n  and for the equality it is sufficient to show that these spaces 
have the same 
dimension. Tensoring with $K$, we have 
$\bigwedge^{2}_{\q}V\otimes_{\q}K=(W_{+}\otimes_{\q}K) \oplus 
(W_{-}\otimes_{\q}K)=W_{+,K}\oplus W_{-,K}$, and looking at the 
decomposition of $V_{K}$ and using the equality 
$\bigwedge^{2}_{\q}V \otimes _{\q} K =\bigwedge^{2}_{K}V_{K}$, we have 
also 
$$\bigwedge^{2}_{\q}V\otimes_{\q}K=\bigwedge^{2}_{K}V_{+} \oplus 
\bigwedge^{2}_{K}V_{-} \oplus \big( V_{+} \otimes_{K} V_{-} \big) .$$
It is easy to show that $\bigwedge^{2}_{K}V_{+} \subset W_{-,K}$, 
$\bigwedge^{2}_{K}V_{-} \subset W_{-,K}$ and $V_{+} \otimes_{K} V_{-} 
\subset W_{+,K}$, and hence we have 
$$
W_{-,K}=\bigwedge^{2}_{K}V_{+} \oplus \bigwedge^{2}_{K}V_{-}, \qquad
W_{+,K}= V_{+} \otimes_{K} V_{-}.
$$
Now, using also the injectivity of the map $i$, we can compute the 
required dimensions:

$${\rm dim}_{K}W_{-,K}=\left( \begin{matrix} n \\ 2 \end{matrix} \right) 
+\left( \begin{matrix} n \\ 2\end{matrix} \right)=2\left( 
\begin{matrix} n \\ 2 \end{matrix} \right), \qquad
{\rm dim}_{K}({\rm Im}(i)\otimes_{\q}K)=2{\rm dim}_K \bigwedge_{K}^2 V =2\left( 
\begin{matrix} n \\ 2 \end{matrix} \right).
\qed $$

\tenv{Abelian varieties of Weil type} 
Let $(X,K,E)$ be a polarized abelian fourfold of Weil type (where 
$K=\Q$ with $\f^2=-d$); we have 
$H^{1}(X,\q)\cong K^{4}$, 
$H^{2}(X,\q)=\bigwedge^{2}_{\q}H^{1}(X,\q)$. 
Now, we study the sub-Hodge structures of the second 
cohomology group of a generic abelian fourfold of Weil type and 
show 
that, if ${\rm discr}(X,K,E)=1$, there is a substructure of weight 2 and 
type (1,4,1).

Let $S=i \bigg( \bigwedge_{K}^2H^1 
(X,\q) \bigg) \subset H^{2}(X, \q )$. From Proposition \ref{naturale} 
the map $i$ is injective, so $S \cong  \bigwedge_{K}^2H^1 
(X,\q)$, and it is often useful to forget the inclusion. Then we have the 
following:

\tenv{\L} \label{282} {\it Let $(X,K,E)$ be an abelian fourfold of 
Weil 
type. Then the subspace $S$ is a 
sub-Hodge structure of $H^2 (X,\q)$ with ${\rm dim}S^{2,0}={\rm dim}S^{0,2}=2, \ 
{\rm dim}S^{1,1}=8$}.

\ts We consider the automorphism $\f_2$ of $H^{2}(X,\q)$;  we have 
$\f_2^2=d^2$, and let $H^2_+$ and $H^2_-$ be its 
(respectively) $+d$ and $-d$ eigenspaces. Using Proposition 
\ref{naturale}, we have $S=H^2_-$, and hence it is a sub-Hodge 
structure of 
$H^2 (X,\q)$. 
Let $H^1(X,\CC)=V_+ \oplus V_-$ be the decomposition into the $i \sqrt 
d$ and 
the $-i \sqrt d$ eigenspaces of $\f$; since $X$ is of Weil type 
$(2,2)$, we have $V_{\pm}=V_{\pm}^{1,0} 
\oplus 
V_{\pm}^{0,1}$, where $V_{\pm}^{i,j}=V_{\pm} \cap H^{i,j}$ and 
${\rm dim}_{\CC}V_{\pm}^{i,j}=2$.
Then we have
$$ \begin{array}{lclcl}
{\rm dim}S^{2,0}&=&{\rm dim}\bigwedge^{2}_{\CC}V_{+}^{1,0}+{\rm dim}\bigwedge^{2}_{\CC}V_{-}^{1,0}&=&2,\\
{\rm dim}S^{0,2}&=&{\rm dim}\bigwedge^{2}_{\CC}V_{+}^{0,1}+{\rm dim}\bigwedge^{2}_{\CC}V_{-}^{0,1}&=&2,\\
{\rm dim}S^{1,1}&=&{\rm dim} V_{+}^{1,0}\otimes 
V_{+}^{0,1}+{\rm dim}V_{-}^{1,0}\otimes 
V_{-}^{0,1}&=&8. 
\end{array}
$$ \qed

\tenv{Hermitian form} We can extend $H$ to 
a $\q$-bilinear form $\H$ on $\bigwedge^{2}_K H^1(X, \q)$ by defining
$$\H (a \w_{K} b, c \w_{K} d) 
:=H(a,c)H(b,d)-H(a,d)H(b,c).$$
 
\n Let $\{v_1, w_1, v_2, w_2\}$ be a $K$-basis of 
$H^1 (X,\q)$ in which the matrix of the Hermitian form 
 is ${\rm diag}(a,1,-1,-1)$. A direct computation shows that in the $K$ 
basis of $\bigwedge_{K}^2H^1(X,\q)$ 
$$a_1=v_1 \w_{K} w_1,\ a_2=v_1 \w_{K} v_2,\ a_3=v_1 \w_{K} w_2,\ 
b_1=v_2 \w_{K}
w_2,\ b_2=w_2 \w_{K} w_1,\ b_3=w_1 \w_{K} v_2,$$
and hence we have $\H \cong {\rm diag} (a,-a,-a,1,-1,-1)$. 

Let now $\gamma :\bigwedge_{K}^4H^1(X,\q)  \longrightarrow  K$ be 
the isomorphism sending $a_1 \wedge_{K} b_1 $ to $1$. Then we have the following

\tenv {\P} {\label{t}} {\it There exists a $\f$-antilinear 
automorphism $t \in {\rm End_{Hod}}(S)$  such that 
$t \circ t=a \cdot {\rm Id}$, and for all 
$v,w \in \bigwedge_{K}^2H^1(X,\q)$ we have $\H (v,w)=\gamma (t(w) 
\wedge v)$  . }

\ts We define a $K$-linear isomorphism  $\rho$
$$
\rho:\bigwedge_{K}^2H^1(X,\q)\longrightarrow
(\bigwedge_{K}^2H^1(X,\q))^*,\qquad
 x\longmapsto\lbrack y\mapsto \gamma
(x\wedge_{K} y)\rbrack 
$$
and a $\f$-antilinear bijection $\tau$
$$
\tau:\bigwedge_{K}^2H^1(X,\q)\longrightarrow
(\bigwedge_{K}^2H^1(X,\q))^*,\qquad
 x\longmapsto\lbrack y\mapsto 
\H (y,x) \rbrack.
$$

\n The isomorphism $t$ is defined as 
$$t:=\rho^{-1} \circ \tau:\bigwedge_{K}^2H^1(X,\q) \longrightarrow
\bigwedge_{K}^2H^1(X,\q).$$
\n Since $\tau (w)= \rho (t(w))$, we have $\H (v,w)=\gamma (t(w) 
\wk
v)$ and the $\f$-antilinearity follows from the $K$-linearity of 
$\rho^{-1}$ and the $\f$-antilinearity of $\tau$.
The representation $h:\CC^* \longmapsto GL(H^1(X,\q))$ defining 
the Hodge structure of $H^1(X,\q)$ gives a representation 
$$
h_2:\CC^* \longmapsto  GL(S),\qquad
  z \longmapsto \lbrack v \wk w \mapsto h(z)v \wk h(z) w \rbrack ,
$$
which defines a Hodge structure on $S$. Since 
$\bigwedge^4_KH^1(X,\q) \cong K$ is one-dimensional, we have that 
$h_4(z)(v \wk w)=h_2(z)v \wk h_2(z)w= |z|^4(v \wk w)$ for all $v,w 
\in S$ and, by the properties of $E$ and $H$, we also have that 
$\H (h_2(z)v,h_2(z)w)=(z \bar{z})^2 \H(v,w)=|z|^4 \H(v,w)$.
These observations show that for all $v,w \in S$
$$\begin{array}{cclcl}
\gamma (t(h_2(z)(v))\wk h_2(Z)w)&=&\H (h_2(z)w,h_2(z)v) 
 &=&(z \bar{z})^2 \H(w,v) \\
 &=&|z|^4 \gamma (t(v) \wk w) 
 &=& \gamma (h_2(z)t(v)\wk h_2(z)w).
\end{array}$$
Hence $t \circ h_2(z)=h_2(z) \circ t$, so $t \in {\rm End_{Hod}}(S)$.
If we write explicitly the action of $t=\rho^{-1} \tau$ on the 
elements of the basis $\{ a_{1}, \dots , b_{3} \}$, then we have 
$$ a_1  \longmapsto aa_1^*  \longmapsto  ab_1,  \qquad b_1  
\longmapsto  
b_1^*   \longmapsto  a_1   $$
\n (and similarly for the other elements of the basis), and hence $t 
\circ t=a \cdot {\rm Id}$ as required. \qed 

\tenv{\C} {\it ${\rm End_{Hod}} (S)\cong {\bf H}$,  where $${\bf H}=\{ \lambda_{1}+\lambda_{2}\f 
+\lambda_{3}t+\lambda_{4} \f \circ t \ ; \ \f^2=-d, \ t^2=a, \ \f 
\circ t=-t \circ \f \}$$ is a quaternion algebra 
over $\q$.
In particular, if $a \notin {\rm Nm}_{K/Q}(K^{*})$, then ${\bf H}$ is a skew 
field and hence $S$ is a simple Hodge structure.}

\ts By Proposition \ref{t} we 
have that $t \in {\rm End_{Hod}}(S)$. Moreover, 
$$h_{2}(\f x \wk y)=h(\f x) 
\wk h(y)=\f h(x) \wk h(y)=\f h_{2} (x \wk y),$$
so $K \hookrightarrow {\rm End_{Hod}}(S)$. Since $t$ is $\f$-antilinear, we 
have ${\bf H} \subset {\rm End_{Hod} }(S)$. 

Let ${\rm MT} \in GL(H^1(X,\q))$ be the Mumford-Tate group of the Hodge 
structure $H^1(X,\q)$
(for the definitions and properties of the Mumford-Tate group, see 
\cite[par. 6.4]{cime}). Then we have
$${\rm End_{Hod}}(S ) \otimes 
\CC = ({\rm End}_{\CC}(S \otimes \CC ))^{{\rm MT}(\CC )}=
{\rm End} \bigg( \bigwedge_{\CC}^2 V_{+} \oplus  
\bigwedge_{\CC}^2 V_{-} \bigg)^{{\rm MT}(\CC)},$$
where $V_{+}$ and $V_{-}$ are the standard and the dual 
representations of ${\rm MT}(\CC) \cong 
SL(4,\CC)$ (the isomorphism ${\rm MT}(\CC)\cong SL(2n, \CC)$ holds for a general polarized Abelian variety of Weil type of dimension $2n$, see \cite{cime}). 
Since $\bigwedge^{4}_{\CC}V_{-}\cong \CC$, we have
$\bigwedge^{2}_{\CC}V_{+}\cong 
\bigwedge^{2}_{\CC}(V_{-})^{*} \cong \bigwedge^{2}_{\CC}V_{-}$,
so $${\rm End_{Hod}}(S) \otimes \CC \cong M_{2}\bigg( {\rm End}\bigg( 
\bigwedge^{2}_{\CC}V_{+} \bigg)^{{\rm MT}(\CC )} \bigg) \cong M_{2}(\CC)$$

\n (the last isomorphism comes from Schur's lemma). Therefore, 
${\rm End_{Hod}} (S)\cong {\bf H}$ because they are of the same dimension.\qed 

\tenv{Hyperbolic lattice} The bilinear form associated to the lattice ``hyperbolic plane" has a matrix ${\rm Hyp}= \left( \begin{matrix} 0&1 \\ 1&0 
\end{matrix} \right)$ (see \cite[p. 14]{bpv}).

\tenv{\T} \label{ker} {\it Let $(X,K,E)$ be an abelian fourfold of 
Weil 
type with discriminant one. Then we have {\rm :}
\begin{enumerate}
\item $S \cong T_{+} \oplus T_{-}$, where $T_{\pm}={\rm Ker} (t \pm {\rm Id})$ are 
sub-Hodge structures.
\item $\varphi :T_{+} {\stackrel{\sim}{\longrightarrow}}  T_{-}$ is 
an isomorphism of Hodge structures.
\item Let $T=T_{+}$. The polarization $\bigwedge^{2}E$ of 
$H^{2}(X,\q)$ induces a polarization $$\bigwedge^2E|_{T \times 
T}=-{1\over
{2d}} \H |_{T \times T}:T \times T \longrightarrow \q ,$$
and we have $\H|_{T \times T}={\rm Hyp} \oplus {\rm Hyp} \oplus \lbrack -2 \rbrack
\oplus \lbrack -2d \rbrack .$
\end {enumerate} }

\ts The first result is obvious, since $t^2={\rm Id}$ by Proposition \ref{t}.
The second follows from the $\f$-antilinearity of $t$: if $v \in T_+$,
we have $t(v)=v$ and 
$t(\f v)=-\f t(v)=-\f v$, so $\f v \in T_-$.
In order to prove the third result, we observe that $\H$ is Hermitian, 
since 
$$\begin{array}{lcl}
\overline{\H(a \wk b,c \wk d)}&=&\overline{H(a,c)H(b,d)}- 
\overline{H(a,d)H(b,c)}\\
&=&H(c,a)H(d,b)-H(d,a)H(c,b)\\
&=&\H (c \wk d, a\wk b).\\
\end{array}$$

\n On $T \times T$ the form $\H$ is symmetric, since
$$\H (v,a) =\gamma (t(a) \wk v)= \gamma (a \wk v)=\gamma (v 
\wk a) =\H
(a,v),$$ 

\n and in particular ${\rm Im}(\H)=0$. By Proposition \ref{naturale}, the 
elements 
$a \wk b$ 
of $S$ can be written as  $(1/2)a \w 
b-(1/2d)\f a \w \f b$, and an easy computation shows that 

$$\begin{array}{lcl}\bigwedge^2E(a \wk b,c\wk d)&=&(1/2) \lbrack 
E(a,c)E(b,d)-E(a,d)E(b,c)\\
& &-(1/d)(E(\f a, c)E(\f b,d)-E(\f 
a,d)E(\f
b,c))\rbrack.\end{array}$$
 Howewer, 
$$\begin{array}{lcl}
\H (a \w b, c\w d)&=& H(a,c)H(b,d)-H(a,d)H(b,c)\\
 &=& (E(\f a, c) + \f E(a,c))(E(\f b, d)+ \f E(b,d))\\
  & &-(E(\f a, d)+ \f E(a,d))(E(\f b, c)+\f E(b,c))\\
 &=& \{ E(\f a, c)E(\f b, d)+ \f^2 E(a,c)E(b,d) \\
 & & -E(\f a,d)E(\f 
b,c)-\f^2E(a,d)E(b,c) \}\\
 & &+\f \{E(\f a,c)E(b,d)+E(a,c)E(\f b,d) \\
  & & -E(\f a,d)E(b,c)-E(a,d)E(\f 
b,c) \} ,
\end{array}$$
and on $T\times T$ (where ${\rm Im}(\H)=0$) we have $\bigwedge^2E|_{T 
\times T}=-(1/2d)
\H |_{T \times T}$.
A direct computation shows that $\H|_{T\times T}\cong 
{\rm diag}(2,-2,2d,-2d,-2,-2d)$.
With an appropriate change of basis, we have $\H|_{T \times T}={\rm Hyp }
\oplus {\rm Hyp} \oplus \lbrack -2 \rbrack
\oplus \lbrack -2d \rbrack $ as required. \qed

\tenv{\C} {\it The Hodge structure $T$ of {\rm 3.8.3} has type 
$(1,4,1)$.}

\ts We have from Lemma \ref{282} that the type of $S$ is $(2,8,2)$.
By Theorem \ref{ker} we have that $S \cong T ^{\oplus 2}$, and hence the 
type 
of $T$ is $(1,4,1)$. \qed

\section{Clifford algebras} 

We showed that there is a polarized sub-Hodge 
structure $T$ of type (1,4,1) with the quadratic 
form $Q \cong {\rm Hyp} \oplus {\rm Hyp} \oplus \lbrack -2 \rbrack \oplus 
\lbrack -2d \rbrack$ contained in the second cohomology group of an 
abelian variety of Weil type with discriminant one. Now, using the Kuga-Satake construction, we construct an abelian 
variety 
 from this Hodge structure and show that this is the variety of 
Weil type we started with.

\tenv{\D} (See \cite[p. 301]{FH}) Let $V$ be a vector space over $\q$ of 
dimension 
$n$ and $\psi$ be
a symmetric, nondegenerate bilinear form. The Clifford algebra 
$C_{n}$ 
is the quotient of the tensor algebra by the two-sided ideal 
$I(\psi)$ 
generated by all elements of the form $v \otimes v-\psi(v,v)$.

\tenv{Notation} We write simply the ``product" $v_1 \cdots v_n$ for the 
class of $v_1 \otimes \dots \otimes v_n$ in $C_n$.

\tenv{\D}  The even subalgebra $C^{+}_{n}$ is the algebra generated 
by all linear combinations of products of an even number of elements 
of $V$.

\tenv{Dimensions} \label{clifdim} We have ${\rm dim}_{\q} C_{n}=2^{n}$ 
and ${\rm dim}_{\q}C^{+}_{n}=2^{n-1}$.

\section{Kuga-Satake varieties}
\n Let $(V,h,\psi)$ be a weight $2$ polarized Hodge structure of type 
$(1,n-2,1)$, and let $\{ g_1, \dots ,g_n \}$ be a basis of $V$ in 
which the symmetric bilinear form $Q=-\psi$ is given by (see Lemma 
\ref{pol}) $$Q=d_{1}X_{1}^{2}+d_2X_2^2-d_3X_3^2- \dots -d_nX_n^2 
\qquad (d_i \in \q_{>0}).$$ 

\tenv{Complex torus} \label{J} Let $J=(1 / \sqrt{d_1d_2})g_1g_2$. Then 
we have $J^2=-1$ 
and the left multiplication by $J$ on $C^+_n$ (the even 
Clifford algebra constructed from $(V,Q)$) defines a complex 
structure 
on $C^+_{n,\RR} :=C^+_n \otimes_{\q} \RR$. Let now 
$C^+_{n,\ZZ}$ be the lattice of linear 
combinations of elements of the basis of $C^+_n$ with 
integer coefficients. Then $$A_0=(C^+_{n,\RR},J)/ 
(C^+_{n,\ZZ})$$ is a complex torus.

\tenv{Polarization}  Let ${\rm Tr}(x)$ be the trace of the map ``right 
multiplication 
on $C^{+}_{n}$ by the element $x \in C^{+}_{n}$''.
 We can define a polarization $E$ on the complex
torus $A_0$ (see \cite{K-S}, \cite{vG}) by setting
$E(v,w):={\rm Tr}(\alpha\iota(v)w)$, where $\iota$ is the canonical 
involution $\iota (g_1^{a_1}\ldots g_{n}^{a_n})=g_{n}^{a_n}\ldots 
g_1^{a_1}$ and $\alpha \in C^{+}_{n}$ is an element such that $\iota 
(\alpha)=-\alpha$ and 
 $E(v,Jv)>0$ for all $v$.
 
\tenv{Kuga-Satake variety} The abelian variety $(A_0,E)$ is the 
Kuga-Satake variety associated to the Hodge structure $(V,h,\psi)$ of 
type 
$(1,n-2,1)$.

\section{Abelian fourfolds}
\tenv{Hodge structures of dimension 6} Let $(V,h,\psi)$ be a weight 
$2$ polarized Hodge structure of type 
$(1,4,1)$, and let 
$\{ f_1, \dots ,f_6 \}$ be a basis of $V$ in which the bilinear form 
$Q=- \psi$ 
is given by the matrix ${\rm Hyp} \oplus {\rm Hyp} \oplus \lbrack l \rbrack  
\oplus \lbrack m \rbrack$ (with $l,m<0$).
In the basis $$e_1=f_1 +f_2,\ e_2=f_1-f_2,\ e_3=f_3+f_4,\ 
e_4=f_3-f_4,\ e_5=f_5,\ e_6=f_6,$$
\n the matrix of $Q$ is diagonal and we have $Q \cong {\rm diag} 
(2,-2,2,-2,l,m)$. We construct the Kuga-Satake variety $(A_0,E)$ 
associated to the Hodge structure, which is an abelian variety of 
dimension ${\rm dim}_{\CC}A_0= {\rm dim}_{\CC} C^+_{6, \RR}=16$. In order to examine
the structure of this variety, the following is very important.

\tenv{\T} \label{gl4} {\it Let $V$ be a rational vector space of 
dimension $6$ equipped with a symmetric bilinear form $Q={\rm Hyp} \oplus 
{\rm Hyp} \oplus \lbrack l \rbrack \oplus 
\lbrack m \rbrack$ $(l,m<0)$, and let $C^+_6$ be the even Clifford 
subalgebra generated by $(V,Q)$. Then we have $C^+_6 \cong gl_4(\q 
(\sqrt{-lm}))$.}

\ts Let $\{ e_{1}, \dots, e_{6} \} $ be the ``diagonal'' basis and 
$C_4$ be the Clifford algebra 
generated by $\{ e_{1}, \dots, e_4 \}$. Let $z=e_{1} \dots e_{6}$.
Then $\q (z)$ is the center of  $C^+_6$ and  
$z^{2}=-16lm$. The map  
$$
\f :(C^+_4 \oplus e_{5} C^-_4)\otimes_{\q} \q (z) \longrightarrow 
C^+_{6}, \qquad
 a \otimes \lambda \longmapsto a \lambda
$$ is an isomorphism. Indeed, it is a surjective homomorphism between 
two 
$\q$-algebras of the same dimension. We have 
$$\begin{array}{lcl}
C_4^+ &\cong & {\rm End}(\bigwedge^{\rm ev}W)\oplus 
{\rm End}(\bigwedge^{\rm odd}W), \\
C^-_4 &\cong &{\rm End}(\bigwedge^{\rm ev}W,\bigwedge^{\rm odd}W)\oplus
{\rm End}(\bigwedge^{\rm odd}W, \bigwedge^{\rm ev}W),
\end{array}
$$  
\n where $W=\langle f_{1},f_{3}\rangle $ (cf. \cite[p. 305]{FH}). 
Therefore we can construct the isomorphism 
$$
\mu: C^+_4 \oplus e_5 C^-_4\longrightarrow  gl_4 (\q), \qquad
 (a,e_5x) \longmapsto \left( \begin{matrix}
a_{ee} & -l x_{oe} \\ x_{eo} & a_{oo} \end{matrix} \right ),
$$
where $a_{ee} \in {\rm End} (\bigwedge^{\rm ev}W)$, $x_{oe}
\in {\rm End} (\bigwedge^{\rm odd}W, \bigwedge^{\rm ev}W)$,
$x_{eo} \in {\rm End}(\bigwedge^{\rm ev}W, \bigwedge^{\rm odd}W)$, 
and 

\n $a_{oo} \in
{\rm End}(\bigwedge^{\rm odd}W)$. So, we have 
$$C^+_6 \cong (C^+_4 \oplus 
e_{5} 
C^-_4)\otimes_{\q} 
\q (z) \cong gl_{4}(\q) \otimes_{\q} \q (z) \cong gl_{4}(\q) 
\otimes_{\q} 
\q(\sqrt{-lm}) \cong gl_{4}(\q(\sqrt{-lm})). \qed $$

\tenv{\C : Poincar\'e's decomposition} From Theorem \ref{gl4} and 
\cite{S} 
we have that, in the general case, ${\rm End}_0(A_0)\cong C^+_6 \cong 
gl_4(\q (\sqrt{-lm}))$, and hence
$A_0 \sim A^4$ (Poincar\'e's theorem), where $A$ is a simple Abelian 
variety with ${\rm End}_0(A) \cong \q (\sqrt{-lm})$.

\tenv{\T} \label{fourfold}{\it The abelian variety $A$ is an abelian 
fourfold of Weil type over $K=\q(\sqrt{-lm})$, and 
there exists a basis in which the Hermitian form $H=E(\f x,y) +\f 
(x,y)$  is diagonal with $a=1$ {\rm (}that is, $H \cong {\rm diag} 
(1,1,-1,-1)${\rm )}.}

\ts Let $\beta=(1 / 4)f_{1}f_{2}f_{3}f_{4}$. We have  
$\beta^{2}=\beta$ and one can verify that the map 

$$
\phi:C_{6}^{+}\longrightarrow C_{6}^{+},\qquad
 x \longmapsto x \cdot \beta
 $$ 
\n has kernel of dimension $24$ over $\q$. Hence the image has 
dimension $8$ and we have that the image of ${\rm Im} \phi \otimes \RR$ in 
the Kuga-Satake variety is isomorphic to $ A$. 

A direct computation shows that ${\rm Im} \phi$ has a basis 
$$\begin{array}{llll}
\epsilon_{1}=f_{2}f_{4},&\epsilon_{2}={1 \over 
2}f_{1}f_{2}f_{3}f_{4},&\epsilon_{3}=f_{2}f_{3}f_{4}f_{5},& 
\epsilon_{4}=f_{1}f_{2}f_{4}f_{5},\\
\delta_{1}=f_{2}f_{4}f_{5}f_{6},&\delta_{2}={1 \over 
2}f_{1}f_{2}f_{3}f_{4}f_{5}f_{6},&\delta_{3}=f_{2}f_{3}f_{4}f_{6},& 
\delta_{4}=f_{1}f_{2}f_{4}f_{6}.
\end{array}$$

\n It is easy to show (using the ``diagonal'' basis) that 
$$\begin{array}{llll}
J\epsilon_{1}=-\epsilon_{2},&J\epsilon_{2}=\epsilon_{1},&J\epsilon_{3}=-\epsilon_{4},&J\epsilon_{4}=\epsilon_{3},\\
J \delta_{1}=-\delta_{2}, & J 
\delta_{2}=\delta_{1},&J\delta_{3}=-\delta_{4},&J\delta_{4}=\delta_{3},
\end{array}$$
\n and 
$$\begin{array}{llll}
z\epsilon_{1}=4 \delta_{1},& z\epsilon_{2}=4 
\delta_{2},&z\epsilon_{3}=4l\delta_{3},&z\epsilon_{4}=4l\delta_{4},\\
z \delta_{1}=-4lm\epsilon_{1},&z\delta_{2}=-4lm\epsilon_{2},&z 
\delta_{3}=-4m\epsilon_{3},&z\delta_{4}=-4m\epsilon_{4}.
\end{array}$$

\n The $i$-eigenspace of the map $J|_{{\rm Im} \phi \otimes \CC}$ is 
spanned 
by  $-i \epsilon_{1}+\epsilon_{2},\ -i 
\epsilon_{3}+\epsilon_{4},\ -i \delta_{1}+\delta_{2},\ -i 
\delta_{3}+\delta_{4},$ and the action of $z$ on this subspace is 
given by the matrix $$\left(\begin{matrix}0&-4lm&0&0 \\ 4&0&0&0 \\ 
0&0&0&-4m 
\\ 0&0&4l&0 \end{matrix} \right).$$

\n The eigenvalues of the matrix are $z$ and $\bar z$ both with 
multiplicity $2$. Hence $A$ is an abelian fourfold of Weil type over
$\q(z) \cong \q (\sqrt{-lm})$.

\n We observe that the every admissible complex structure $J'$ can be 
written as $J'=a J \io (a)$, 
where $a \in Spin(Q) 
:= \{x \in 
C^+ \ ; \ x \io (x)=1, \ xV\io (x)\subset V \}$
($J'=e'_1e'_3$ with ${e'_i}^{2}>0$; then  $e'_i$ is obtained from 
$e_i$ by the action of $SO(Q)$ and $Spin(Q)$ is a $2:1$ cover of 
$SO(Q)$ ).
Since $$\begin{array}{cclcl}
E(x,J'x) &=& {\rm }tr(\alpha \io (x)J')
 &=& {\rm }tr(\alpha \io (x)a J \io (a) x) \\
 &=& {\rm }tr(\alpha \io(\io (a) x) J (\io (a)x))
 &=& {\rm }tr(\alpha \io (y)Jy),
\end{array}$$
\n we have that the choice of the element $\alpha \in C^{+}_{6}$ in 
the 
definition of the polarization does not depend on the choice of the 
complex structure ($E(x,Jx)>0$ for all admissible $J$).

A direct computation shows that $\alpha =-f_{1}f_{3}$ satisfies 
the ``positivity condition'' of $E$ and that the matrix of the 
polarization in the basis $\{ \ep_{i}, \delta_{j} \}$ is 

$$E= \left( \begin{matrix} M& & & \\ &-2lM& & \\ & & lmM& \\ & & & 
-2mM \end{matrix} \right), \qquad \qquad {\rm where} \ \ M= \left( 
\begin{matrix} 0& -64 \\ 64 &0 \end{matrix} \right).$$


\n On the $K$-basis $\{ \ep_{1}, \dots , \ep _{4} \}$ 
we 
have that $H(x,y)=E(\f x,y)+\f E(x,y)$ has the matrix representation
$$H=\left( \begin{matrix}
\f M & \\ &-2 \f M  \end{matrix}\right),$$ and this matrix can be 
diagonalized as $H \cong {\rm diag}(1,1,-1,-1)$. \qed

We can now prove the following
\tenv{\T} {\it The abelian fourfold $A$ occurring in the decomposition of the Kuga-Satake variety is isogenous to the Abelian fourfold $X$ we started with.}

\ts The abelian fourfolds $A$ and $X$ are both of Weil type with discriminant equal to one, so we have only to show that they have the same complex structure. We consider the Hodge substructure $T\subset H^2(X,\QQ)$ defined in Theorem \ref{ker}. We also have that $T\subset H^2(A,\QQ)$ (see \cite{cime}). Let ${\rm MT}(X)$ be the Mumford-Tate group of the abelian variety $X$, that is the Mumford -Tate group of the Hodge structure $H^1(X,\QQ)$. The subspace $T$ is obviously a subrepresentation of ${\rm MT}(X)$ and its Hodge structure is given by a representation $h_+: \CC^* \longrightarrow  {\rm MT}(T)(\RR)\subset GL\left( T_{\RR}\right)\subset GL(H^2(X,\RR))$. Now, we consider the map $h_X:\CC^*\longrightarrow {\rm MT}(X)(\RR) \subset GL(H^1(X,\RR))$ which gives the Hodge structure on $H^1(X,\QQ)$; we observe that the map ${\rm MT}(X)(\RR)\longrightarrow {\rm MT}(T)(\RR)$ is a 2:1 cover (over $\CC$, and with $A$ and $X$ general, we have that it is the map $SL(4)\longrightarrow SO(6)\cong SL(4)/\langle \pm {\rm Id} \rangle)$ given by the action of ${\rm MT}(X)(\RR)$ on $\bigwedge^2H^1(X,\RR)$. We have then the diagram
$$\begin{array}{cccll}
h_+: & \CC^*& \longrightarrow &{\rm MT}(X)(\RR)/\langle {\pm \rm Id} \rangle={\rm MT}(T)(\RR)& \subset GL(H^2(X,\RR)) \\ & & & \qquad \uparrow & \\
h_X: & \CC^*& \longrightarrow &{\rm MT}(X)(\RR)& \subset GL(H^1(X,\RR)).
\end{array}$$
The complex structure $J_X$ on $X$ is given by $h_X(i)$, which lies over $h_+(i)$, so we have two possible choices for $J_X$, $J$ and $-J$. We repeat now the same argument using the inclusion $T \subset H^2(A,\QQ)$; from ${\rm MT}(X)={\rm MT}(A)$ (see \cite{cime}) we obtain the diagram 
$$\begin{array}{cccll}
h_+: & \CC^*& \longrightarrow &{\rm MT}(X)(\RR)/\langle {\pm \rm Id} \rangle={\rm MT}(T)(\RR)& \subset GL(H^2(A,\RR)) \\ & & & \qquad \uparrow & \\
h_A: & \CC^*& \longrightarrow &{\rm MT}(A)(\RR)={\rm MT}(X)(\RR)& \subset GL(H^1(A,\RR))
\end{array}$$
where $h_A$ gives the Hodge structure on $H^1(A,\QQ)$. Then, for the complex structure on $A$, $J_A=h_A(i)$, we have the same two choices $J$ and $-J$ as $X$ since also $h_A(i)$ lies over $h_+(i)$. As $A$ and $X$ are of Weil type, we can identify the $K$-vector spaces $H^1(A,\QQ)$ and $H^1(X,\QQ)$, and thus $J_A=\pm J_X$. Since the polarization on an abelian fourfold of Weil type with discriminant equal to one is unique (see \cite{cime}), we have $E_X=E_A$. From the ``positivity condition'' of a polarization we conclude that $A$ and $X$ have the same complex structure. \qed

\section{Higher dimensions}

Now, we generalize the result of Theorem 
\ref{fourfold} and show that there exist other cases in which the 
Kuga-Satake construction 
gives abelian varieties of Weil type with discriminant one.

Let $(V,h,Q)$ be a rational polarized weight $2$ Hodge structure of 
${\rm dim}_{\q}V=n=2m$, type $(1,n-2,1)$, and let $Q=-\psi$ given by 
$Q=d_{1}X_{1}^{2}+d_2X_2^2-d_3X_3^2- \dots -d_nX_n^2$ ($d_i \in 
\q_{>0}$) in a basis $\{g_{1}, \dots, g_{n} \}$.
Let $-d:=(-1)^{m}d_{1} \dots  d_{n}$ and $C:=C_n$ be the Clifford 
algebra associated to $(V,Q)$. We consider the case 
$n \equiv 2 \pmod 4$ equivalent to (see \cite[Theorem 7.7]{vG}) $C^{+} \cong M_{2^{m-1}}(\q (\sqrt{-d}))$  and construct the 
Kuga-Satake variety $A_0$ associated to such a Hodge structure.
Since ${\rm End}(A_0) \cong C^{+} \cong M_{2^{m-1}}(\q(\sqrt{-d}))$, we 
have from Poincare's theorem that 
$A_0 \sim A^{2^{m-1}}$, where $A$ is an abelian variety with ${\rm End}(A) 
\cong 
Q(\sqrt{-d})$. Obviously, we have that 
${\rm dim}_{\RR}C^+_{\RR}={\rm dim}_{\q}C^+=2^{2m-1}$, and therefore 
${\rm dim}(A_0)={\rm dim}_{\CC}C^+_{\RR}=2^{2m-2}$ and 
${\rm dim}_{\CC}(A)=2^{m-1}$.

\tenv{\T} {\it Let $(V,h,\psi)$ be a rational polarized weight $2$ 
Hodge structure with dimension 

\n ${\rm dim}_{\q}V=n=2m$, type $(1,n-2,1)$
 and let $-d:=(-1)^{m}d_{1} \dots  d_{n}$. Then, if $n \equiv 2 \pmod 4$, $A$ is an abelian variety of Weil type over $K=\q(\sqrt {-d})$.}

\ts
First, we consider the case in which $J=(1/ 
\sqrt{d_{1}d_{2}})g_{1}g_{2}$ is the complex structure on 
$C^+_{\RR}$. 
Let $z=g_{1}  \dots  g_{n}$. Then the element $z$ is in the center of 
$C^{+}$ and we have 
 $$\begin{array}{lclcl}
z^{2}&=&g_{1} \dots  
g_{n} g_{1}  \dots  g_{n}&=&(-1)^{(n-1)+(n-2) \dots 
1}g_{1}^{2} \dots g_{n}^{2}\\
&=&(-1)^{m(2m-1)}d_{1} \dots d_{n}
&=&-d 
\end{array}$$
(since $2m-1$ is odd). Let $x \in C^{+}$ a generic element of the even subalgebra. We then  
define a map ``left multiplication for $x$'' by setting 

$$l_x:C^{+}_{\RR}\longrightarrow  C^{+}_{\RR},  \qquad
 v  \longmapsto  xv. $$

\n The maps $l_J$ and $l_z$ are injective (e.g.,  
$l_z(w_{1})=l_z(w_{2}) \Rightarrow  
l_z^2(w_{1})=l_z^2(zw_{2})\Rightarrow w_{1}=w_{2}$) and both have 
complex eigenvalues. Therefore
we have that, for all $v \in C^+_{\RR}$,  $l_J(v) \notin \langle v \rangle$ and 
$l(z)v \notin \langle v \rangle $.
Let ${\bf H}=\q(J,z)$; ${\bf H}$ is an extension of degree 4 over $\q$
and there are the obvious isomorphisms $H^{1}(A_{0}, {\RR})\cong 
C^{+}_{\RR} \cong {\bf H}^{2m-3}$.
\n Using these isomorphisms, we can construct a basis of
$H^{1}(A_{0}, \RR )$ in the following way: let $\widetilde g_1=1 \in 
C^+$. Using the ``left multiplications",  we obtain $J$, $z$ and $Jz$ 
(these elements are independent). 

Now we choose an element 
$\widetilde g_2$ not contained in the span of the previous elements 
and  
continue the construction. In this  way we find the basis
$$\{\widetilde g_1=1,J, z, Jz, \dots ,\widetilde g_t, J \widetilde 
g_t, z\widetilde g_t, Jz\widetilde g_t \} \qquad (t=2^{2m-1} / 
4=2^{2m-3}).$$

\n It is easy to show that 
$$\{ i\widetilde g_1+J\widetilde g_1, iz\widetilde g_1+Jz\widetilde 
g_1, \dots ,i\widetilde g_t+J\widetilde g_t, iz\widetilde g_t+
Jz\widetilde g_t \}$$

\n is a basis of $H^{1,0}(A_0)$ (the $i$-eigenspace of $J$), and we 
can compute directly the action of $z$ on this space:

\n $$\begin{array}{lcl}
z(ig_{j}+Jg_{j})&=&izg_{j}+Jzg_{j}\\
z(izg_{j}+Jzg_{j})&=&d(ig_{j}+Jg_{j}) \\
\end{array}, 	\quad
{\rm and \ hence} \quad
z \cong \left( \begin{matrix} {\begin{matrix} 0&-d \\ 1&0 
\end{matrix}} & & \\  &\ddots & 
\\ & & {\begin{matrix}0&-d \\ 1&0 \end{matrix}}\end{matrix} 
\right).$$

\n Since every  $$\left(\begin{matrix}0&-d \\ 1&0 
\end{matrix} \right) \mbox{can be diagonalized as }
\left( \begin{matrix} \sqrt {-d} &0 \\ 0&- \sqrt {-d} 
 \end{matrix}\right),$$ 
we see that $A$ is of Weil type.

Now, we show that this result holds for all complex structures. 
Let $C^{+}_{\CC}=V_{+} \oplus V_{-}$ be the decomposition of 
$C^{+}_{\CC}$ in the $\left( \sqrt{-d}\right)$- and  $\left( -\sqrt{-d}\right)$-eigenspaces of 
$z$. We observe that $J$ commutes with $z$ ($z$ is contained in the 
center), so $l_J$ respects this decomposition and it has $t$ 
eigenvalues $i$ and $t$ eigenvalues $\bar i$ on each component 
(indeed, $A$ is of Weil type and the condition ``$l_z$ has type 
$(t,t)$ 
on $H^{1,0}(X)$'' is obviously equivalent to the condition ``$l_{J}$ 
has type $(t,t)$ on $V_{+}$''). Moreover, every complex structure can 
be written as $gJg^{-1}$ with $g \in Spin(Q)$; we have $gz=zg$, so 
$l_g$ respects the decomposition of $ C^{+}_{\CC}$, and therefore 
$l_{gJg^{-1}}$ respects this decomposition. Since the matrices of  
$l_J$ and $l_{gJg^{-1}}$ have the same eigenvalues, $l_{gJg^{-1}}$ 
has type $(t,t)$ and $A$ is of Weil type for every complex structure.
\qed

\subsection{Discriminant.}

\n Let $E(x,y)={\rm Tr} (\alpha \iota (x)y)$ be the polarization of $A$. We want to show that ${\rm discr}(A,\q(\sqrt {-d}),E)=1$.
 
\tenv {\T} {\it There exists a $\q(\sqrt {-d})$ basis of $H_1(A,\q)$ 
in which  
$$H \cong {\rm diag} (1, \dots , 1, -1, \dots , -1).$$}
\ts
Using the decomposition $$C^+=(C^{+}_{n-2} \oplus g_{n-1} 
C^{+}_{n-2})\otimes_{\q} \q (z),$$
it is possible to construct a  $\q(z)$-basis 
$\{ e_1, \dots, e_h \} \quad ( {\rm where} \quad h=2^{n-2})$ of 
$H_1(A,\q)$
in such a way that  $g_{n}$ does not appear in the elements $e_i$ (as 
complex structure we can use, for example, the element $J$ defined in 
\ref{J} ).
Using the Frobenius theorem, we can change the basis $\{ e_i \}$ to a 
basis $\{ \widetilde e_i \}$ in which the matrix of the polarization 
has the form 
$$E\sim \left( \begin{matrix} { \begin{matrix} 0&a_{1} \\ -a_{1}&0 
\end{matrix}} & & \\  &\ddots 
& \\ & & {\begin{matrix} 0&a_{h} \\ -a_{h}&0 \end{matrix} 
} \end{matrix}\right) \qquad (a_1 \in \q).$$
 
\n Now, we decompose the even Clifford subalgebra as $C^{+}=\q \oplus 
C^{+}_{0}$; from  \cite {vG} we have that ${\rm Ker(Tr)}= C^{+}_{0}$. Then, 
we can choose the element $\alpha \in 
C^{+}_{n}$ of the polarization without terms containing $g_{n}$.
Indeed, $\iota (x)y$ does not contain $g_{n}$, and if 
$\alpha=\alpha_{1}+g_{n}\alpha_{2}$, we have 
$${\rm Tr}(\alpha \iota (x)y)={\rm Tr}(\alpha_{1}\iota 
(x)y)+{\rm Tr}(g_{n}\alpha_{2}\iota (x)y)={\rm Tr}(\alpha_{1}\iota (x)y)+0$$

\n ( $g_{n}\alpha_{2}\iota (x)y$ contains necessarily $g_{n}$, so it 
cannot be a coefficient), so the term $g_{n}\alpha_{2}$ is useless.

\n By definition, $H(x,y)=E(zx,y)+\sqrt{-d}E(x,y)$, 
and we have $E(z\widetilde{e_i}, \widetilde{e_j})={\rm Tr}(a\iota 
(z\widetilde{e_i})z\widetilde{e_j})=0$ (the argument contains 
$g_{n}$). Then

$$H\cong \left( \begin{matrix} \begin{matrix} 0&\sqrt{-d}a_{1}  \\ 
-\sqrt{-d}a_{1}&0 \end{matrix} & & \\  &\ddots 
& \\ & & \begin{matrix} 0&\sqrt{-d}a_{h} \\ -\sqrt{-d}a_{h}&0 
\end{matrix} \end{matrix} \right).$$

On a $\q(\sqrt{-d})$-basis, this matrix can be transformed as

$$H \cong \left( \begin{matrix}  \begin{matrix} 0&1 \\ 1&0 
\end{matrix} & & \\  &\ddots & 
\\ & & \begin{matrix}0&1 \\ 1&0 \end{matrix} \end{matrix} 
\right)$$

\n and we have $H\cong {\rm diag}(1,\dots 
,1,-1,\dots ,-1)$ as required. \qed

\vskip1truecm

{\footnotesize \sc
Lombardo Giuseppe

Dipartimento di Matematica

Universit\'a di Torino 

Via Carlo Alberto 10

10123 TORINO

ITALIA

\smallskip

{\it E-mail address}{\rm : lombardo@dm.unito.it}}
\end{document}